\documentclass[11pt]{article}
\usepackage{amsmath,amssymb,latexsym, amsfonts, amscd, amsthm}
\usepackage[small,nohug,heads=littlevee]{diagrams}
\diagramstyle[labelstyle=\scriptstyle]

\theoremstyle{plain}

\newtheorem{theorem}{Theorem}[section]

\newtheorem{proposition}[theorem]{Proposition}

\newtheorem{lemma}[theorem]{Lemma}

\theoremstyle{definition}
\newtheorem{definition}[theorem]{Definition}

\newcommand{\mat}[4]{\left( \begin{array}{cc} {#1} & {#2} \\ {#3} & {#4}
\end{array} \right)}

\def\bdf{\begin{defn}}
\def\edf{\end{defn}}

\usepackage{color}

\begin{document}

\title{On the Local Langlands Correspondences of DeBacker--Reeder and Reeder for $GL(\ell,F)$, where $\ell$ is prime}

\author{Moshe Adrian}

\maketitle

\begin{abstract}
We prove that the conjectural depth-zero local Langlands correspondence of DeBacker--Reeder agrees with the known depth-zero local Langlands correspondence for the group $GL(\ell,F)$, where $\ell$ is prime and $F$ is a nonarchimedean local field of characteristic 0.  We also prove that if one assumes a certain compatibility condition between Adler's and Howe's construction of supercuspidal representations, then the conjectural positive depth local Langlands correspondence of Reeder also agrees with the known positive depth local Langlands correspondence for $GL(\ell,F)$.
\end{abstract}

\section{Introduction}

Let $F$ be a nonarchimedean local field of characteristic zero.  Let $\mathbf{G}$ be a connected reductive group defined over $F$.  The local Langlands correspondence asserts that there is a finite to one map from the set of admissible representations of $\mathbf{G}(F)$ to the set of Langlands parameters of $\mathbf{G}(F)$, satisfying various conditions.  Until recently, this has only been proven for special cases of groups such as $GL(n,F)$, $Sp(4,F)$, and $U(3)$.
The local Langlands correspondence for $GL(n,F)$ has been proven independently by Harris/Taylor and Henniart.  More recently, DeBacker--Reeder and Reeder, in \cite{debackerreeder} and \cite{reeder}, have described conjectural local Langlands correspondences for a more general class of groups, and certain classes of Langlands parameters.  Despite the fact that their correspondences satisfy several requirements that the Langlands correspondence should have, their correspondences are still conjectural.  One would therefore like to know whether they agree at least with the proven correspondences in the known cases.  We prove that the correspondence of \cite{debackerreeder} agrees with the known correspondence for $GL(\ell,F)$, and that the correspodnence of \cite{reeder} agrees with the known correspondence for $GL(\ell,F)$ if one assumes a certain compatibility condition which we describe later.

For $GL(n,F)$, the constructions of Harris/Taylor, Henniart, and DeBacker--Reeder (and Reeder) use different methods.  We first recall the classical construction of the tame local Langlands correspondence for $GL(\ell,F)$ as in \cite{moy}.  We note that a tame local Langlands correspondence for $GL(n,F)$ was conjectured by Moy in \cite{moy} for general $n$. Because of recent work of Bushnell and Henniart (see \cite{bushnellhenniart1}), the correspondence of \cite{moy} is indeed correct for $GL(\ell,F)$, $\ell$ a prime.

\begin{definition}
Let $E/F$ be an extension of degree $\ell$, $\ell$ relatively prime to the residual characteristic of $F$, and let $\chi$ be a character of $E^*$.  The pair $(E/F, \chi)$ is called \emph{admissible} if $\chi$ does not factor through the norm from a proper subfield of $E$ containing $F$.
\end{definition}

We write $\mathbb{P}_{\ell}(F)$ for the set of $F$-isomorphism classes of admissible pairs $(E/F, \chi)$ where $E/F$ is a degree $\ell$ extension (for more information about admissible pairs, see \cite{moy}).  Let $\mathbb{A}_{\ell}^0(F)$ denote the set of supercuspidal representations of $GL(\ell,F)$. Howe constructs a map (see \cite{howe}) $$\mathbb{P}_{\ell}(F) \rightarrow \mathbb{A}_{\ell}^0(F)$$ $$(E/F, \chi) \mapsto \pi_{\chi}$$  This map is a bijection (see \cite{moy}).  Let $\mathbb{G}_{\ell}^0(F)$ denote the set of irreducible $\ell$-dimensional representations of $W_F$, where $W_F$ is the Weil group of $F$.  We then have a bijection (see \cite{moy}) $$\mathbb{P}_{\ell}(F) \rightarrow \mathbb{G}_{\ell}^0(F)$$ $$(E/F, \chi) \mapsto Ind_{W_E}^{W_F}(\chi) =: \phi(\chi)$$  The local Langlands correspondence is then given by $$\phi(\chi) \mapsto \pi_{\chi \Delta_{\chi}}$$ for some subtle finite order character $\Delta_{\chi}$ of $E^*$ (see \cite{bushnellhenniart1}.  In the case of depth-zero supercuspidal representations, there is only one extension $E/F$ to deal with, namely, the unramified extension of $F$ of degree $\ell$.

On the other hand, the constructions of DeBacker--Reeder (see \cite{debackerreeder}) and Reeder (see \cite{reeder}) extensively use Bruhat-Tits theory.  To a certain class of Langlands parameters for an unramified connected reductive group $\mathbf{G}$, they associate a character of a torus, to which they attach a collection of supercuspidal representations on the pure inner forms of $\mathbf{G}(F)$, a conjectural $L$-packet.  They are also able to isolate the part of this packet corresponding to a particular pure inner form, and prove that their correspondences satisfy various natural conditions, such as stability.

Specifically, we will prove the following.  Let $E/F$ be the unramified degree $\ell$ extension, $\ell$ a prime. To any tame, regular, semisimple, elliptic, Langlands parameter (TRSELP) for $GL(\ell,F)$, we show that the theory of DeBacker--Reeder attaches the character $\chi \Delta_{\chi}$ of $E^*$, to which is attached the representation $\pi_{\chi \Delta_{\chi}}$.  This will prove that their correspondence agrees with the correspondence of \cite{moy} for $GL(\ell,F)$.

We then prove the same for Reeder's construction, if one assumes a certain compatibility condition, which we describe now.  Reeder's construction in \cite{reeder} begins by canonically attaching a certain admissible pair $(L/F, \Omega)$ to a Langlands parameter for $GL(\ell,F)$.  His construction then inputs this admissible pair into the theory of \cite{adler} in order to construct a supercuspidal representation $\pi(L, \Omega)$ of $GL(\ell,F)$.  The compatibility condition that we will need to assume is that $\pi(L, \Omega)$ is the same supercuspidal representation that is attached to $(L/F, \Omega)$ via the Howe construction in \cite{howe}.  We remark that this compatibility condition does not seem to be known to the experts.

We wish to make an interesting remark regarding the correspondences of Moy and \\ DeBacker--Reeder (and Reeder, assuming the above compatibility).  Although the correspondence of Moy agrees with DeBacker--Reeder (and Reeder, assuming the above compatibility), there are some important details that are different.  In particular, the passage from a Langlands parameter to a character of a torus has a subtle but important difference.  To illustrate this difference, we rewrite both correspondences to include their factorization through characters of elliptic tori as
$$ \{ \mathrm{Langlands \ Parameters \ from} \ \cite{debackerreeder} \ \mathrm{or} \ \cite{reeder} \ \mathrm{for} \ GL(\ell,F)  \} \rightarrow \mathbb{P}_{\ell}(F) \rightarrow \mathbb{A}_{\ell}^0(F)$$
Then, the correspondence of Moy is given by $$\phi(\chi) = Ind_{W_E}^{W_F}(\chi) \mapsto (E/F, \chi) \mapsto \pi_{\chi \Delta_{\chi}}$$ whereas the correspondences of DeBacker--Reeder (and Reeder, assuming the compatibility) are given by $$\phi(\chi) = Ind_{W_E}^{W_F}(\chi) \mapsto (E/F, \chi \Delta_{\chi}) \mapsto \pi_{\chi \Delta_{\chi}}.$$

We now briefly present an outline of the paper.  In section \ref{notation}, we introduce some notation that we will need throughout.  In section \ref{preliminaries}, we briefly recall some of the key components to the construction of DeBacker--Reeder from \cite{debackerreeder}.  In section \ref{existingdescription}, we recall the tame local Langlands correspondence for $GL(\ell,F)$ as explained in \cite{moy}.  In sections \ref{glellf} and \ref{glellfcharactertorus}, we work out the theory of \cite{debackerreeder} for $GL(\ell,F)$, and we show that the correspondence of DeBacker--Reeder and Moy agree for $GL(\ell,F)$.  Finally, in section \ref{positivedepthreederchapter}, we work out the theory of \cite{reeder} for $GL(\ell,F)$, where $\ell$ is prime, and we show that under the compatibility condition, the correspondence of Reeder and Moy agree for $GL(\ell,F)$.

\section{Notation}\label{notation}
Let $F$ denote a nonarchimedean local field of characteristic zero.  We let $\mathfrak{o}_F$ denote the ring of integers of $F$, $\mathfrak{p}_F$ its maximal ideal, $\mathfrak{f}$ the residue field of $F$, $q$ the order of $\mathfrak{f}$, and $p$ the characteristic of $\mathfrak{f}$.  Let $\mathfrak{f_m}$ denote the degree $m$ extension of $\mathfrak{f}$.  We let $\varpi$ denote a uniformizer of $F$.  Let $F^u$ denote the maximal unramified extension of $F$.  We have the canonical projection $$\Pi :  \mathfrak{o}_F^* \rightarrow \mathfrak{o}_F^* / (1 + \mathfrak{p}_F) \cong \mathfrak{f}^*$$
We denote by $W_F$ the Weil group of $F$, $I_F$ the inertia subgroup of $W_F$, $I_F^+$ the wild inertia subgroup of $W_F$, and $W_F^{ab}$ the abelianization of $W_F$.  We denote by $W_F'$ the Weil-Deligne group, we set $W_t := W_F / I_F^+$, and we set $I_t := I_F / I_F^+$.  We fix an element $\Phi$ $\in Gal(\overline{F} / F)$ whose inverse induces the map $x \mapsto x^q$ on $\mathfrak{F} := \overline{\mathfrak{f}}$, and if $E/F$ is the unramified extension of degree $\ell$, we fix an element $\Phi_E$ $\in Gal(\overline{E} / E)$ whose inverse induces the map $x \mapsto x^{q^{\ell}}$ on $\mathfrak{F} := \overline{\mathfrak{f}}$.  Let $\mathbf{G}$ be an unramified connected reductive group over $F$, and set $G = \mathbf{G}(F^u)$.  We fix $\mathbf{T} \subset \mathbf{G}$, an $F^u$-split maximal torus which is defined over $F$ and maximally $F$-split, and set $T = \mathbf{T}(F^u)$.  We write $X := X_*(\mathbf{T})$, $W_o$ for the finite Weyl group $N_G(T) / T$, and set $N := N_G(T)$.  Recall that the extended affine Weyl group is defined by $W := X \rtimes W_o$, and that the affine Weyl group is defined by $W^o := \Psi \rtimes W_o$, where $\Psi$ is the coroot lattice in $X$.  We let $\mathcal{A} := \mathcal{A}(T)$ be the apartment of $T$.  We denote by $\theta$ the automorphism of $X, W$ induced by $\Phi$.  If $E/F$ is a finite Galois extension, then we denote by $\aleph_{E/F}$ the local class field theory character of $F^*$ with respect to the extension $E/F$.  If $\chi \in \widehat{E^*}$ satisfies $\chi|_{1 + \mathfrak{p}_E} \equiv 1$, then $\chi|_{\mathfrak{o}_E^*}$ factors to a character, denoted $\chi_o$, of the multiplicative group of the residue field of $E$, given by $\chi_o(x) := \chi(u)$ for any $u \in \mathfrak{o}_E^*$ such that $\Pi(u) = x$.  If $E/F$ is the degree $\ell$ unramified extension, where $\ell$ is prime, we once and for all fix a generator $\xi$ of $Gal(E/F)$.  We also fix a generator of $Gal(\mathfrak{f}_{\ell}/\mathfrak{f})$, which, abusing notation, we also denote by $\xi$.  If $\chi$ is a character of $E^*$ or $\mathfrak{f}_{\ell}^*$, we let $\chi^{\xi}$ denote the character given by $\chi^{\xi}(x) := \chi(\xi(x))$.  If $L / K$ is a Galois quadratic extension, we let the map $x \mapsto \overline{x}$ denote the nontrivial Galois automorphism of $L/K$.  If $A$ is a group and $B$ is a normal subgroup of $A$, we denote the image of $a \in A$ in $A / B$ by $[a]$.  If $\phi : C \rightarrow D$ is a group homomorphism and $\phi$ is trivial on a normal subgroup $M \lhd C$, then we will abuse notation and write $\phi|_{C / M}$ for the factorization of $\phi$ to a map $C / M \rightarrow D$.  For example, the Langlands parameters in \cite{debackerreeder} are trivial on the wild inertia subgroup $I_F^+$ of the inertia group $I_F$.  Therefore, if $\phi$ is such a Langlands parameter and $I_t := I_F / I_F^+$, we will write $\phi|_{I_t}$ to denote the factorization of $\phi|_{I_F}$ to the quotient $I_t$.

\section{Review of Construction of DeBacker and Reeder}\label{preliminaries}

We first review some of the basic theory from \cite{debackerreeder}.  We first fix a pinning $(\hat{T}, \hat{B}, \{x_{\alpha} \})$ for the dual group $\hat{G}$.  The operator $\hat{\theta}$ dual to $\theta$ extends to an automorphism of $\hat{T}$.  There is a unique extension of $\hat{\theta}$ to an automorphism of $\hat{G}$, satisfying $\hat{\theta}(x_{\alpha}) = x_{\theta \cdot \alpha}$ (see \cite[section 3.2]{debackerreeder}).  Following \cite{debackerreeder}, we may form the semidirect product ${}^L G := \ <\hat{\theta}> \ltimes \hat{G}$.

\begin{definition}
Let $W_F'$ denote the Weil-Deligne group.  A Langlands parameter $\phi : W_F' \rightarrow {}^L G$ is called a \emph{tame regular semisimple elliptic Langlands parameter} (abbreviated TRSELP) if

(1) $\phi$ is trivial on $I_F^+$,

(2) The centralizer of $\phi(I_F)$ in $\hat{G}$ is a torus.

(3) $C_{\hat{G}}(\phi)^o = (\hat{Z}^{\hat{\theta}})^o$, where $\hat{Z}$ denotes the center of $\hat{G}$.
\end{definition}

Condition (2) forces $\phi$ to be trivial on $SL(2,\mathbb{C})$.  Let $\hat{N} = N_{\hat{G}}(\hat{T})$.  After conjugating by $\hat{G}$, we may assume that $\phi(I_F) \subset \hat{T}$ and $\phi(\Phi) = \hat{\theta} f$, where $f \in \hat{N}$.  Let $\hat{w}$ be the image of $f$ in $\hat{W}_o$, and let $w$ be the element of $W_o$ corresponding to $\hat{w}$.

Let $\phi$ be a TRSELP with associated $w$ and set $\sigma = w \theta$. $\sigma$ is an automorphism of $X$.  Let $\hat{\sigma}$ be the automorphism dual to $\sigma$, and let $n$ be the order of $\sigma$.  We set $\hat{G}_{ab} := \hat{G} / \hat{G}'$, where $\hat{G}'$ denotes the derived group of $\hat{G}$. Let ${}^L T_{\sigma} := \langle \hat{\sigma} \rangle \ltimes \hat{T}$.  Associated to $\phi$, DeBacker--Reeder (see \cite[Chapter 4]{debackerreeder}) define a $\hat{T}$-conjugacy class of Langlands parameters

\begin{equation}
\phi_T : W_F \rightarrow {}^L T_{\sigma} \ \label{phiT}
\end{equation}

\noindent as follows.  Set $\phi_T := \phi$ on $I_F$, and $\phi_T(\Phi) := \hat{\sigma} \ltimes \tau$ where $\tau \in \hat{T}$ is any element whose class in $\hat{T} / (1 - \hat{\sigma}) \hat{T}$ corresponds to the image of $f$ in $\hat{G}_{ab} / (1 - \hat{\theta}) \hat{G}_{ab}$ under the bijection

\begin{equation}
\hat{T} / (1 - \hat{\sigma}) \hat{T} \stackrel{\sim}{\rightarrow} \hat{G}_{ab} / (1 - \hat{\theta}) \hat{G}_{ab} \ \label{bijectionfortau}
\end{equation}

In \cite[Chapter 4]{debackerreeder}, DeBacker and Reeder construct a canonical bijection between $\hat{T}$-conjugacy classes of admissible homomorphisms $\phi : W_t \rightarrow {}^L T_{\sigma}$ and depth-zero characers of $T^{\Phi_{\sigma}}$ where $\Phi_{\sigma} := \sigma \otimes \Phi^{-1}$.  We briefly summarize this construction. Let $\mathbb{T} := X \otimes \mathfrak{F}^*$.  Given automorphisms $\alpha, \beta$ of abelian groups $A,B$, respectively, let $Hom_{\alpha, \beta}(A,B)$ denote the set of homomorphisms $f : A \rightarrow B$ such that $f \circ \alpha = \beta \circ f$.  The twisted norm map $$N_{\sigma} : \mathbb{T}^{\Phi_{\sigma}^n} \rightarrow \mathbb{T}^{\Phi_{\sigma}}$$ given by $N_{\sigma}(t) = t \Phi_{\sigma}(t) \Phi_{\sigma}^2(t) \cdots \Phi_{\sigma}^{n-1}(t)$ induces isomorphisms $$\mathrm{Hom}(\mathbb{T}^{\Phi_{\sigma}}, \mathbb{C}^*) \stackrel{\sim}{\rightarrow}  \mathrm{Hom}_{\Phi_{\sigma}, Id}(\mathbb{T}^{\Phi_{\sigma}^n}, \mathbb{C}^*) \stackrel{\sim}{\rightarrow} \mathrm{Hom}_{\Phi_{\sigma}, Id}(X \otimes \mathfrak{f}_n^*, \mathbb{C}^*)$$
Moreover, the map $s \mapsto \chi_s$ gives an isomorphism $$\mathrm{Hom}_{\Phi, \hat{\sigma}}(\mathfrak{f}_n^*, \hat{T}) \stackrel{\sim}{\rightarrow} \mathrm{Hom}_{\Phi_{\sigma}, Id}(X \otimes \mathfrak{f}_n^*, \mathbb{C}^*)$$ where $\chi_s(\lambda \otimes a) := \lambda(s(a))$.  The canonical projection $I_t \rightarrow \mathfrak{f}_m^*$ induces an isomorphism as $\Phi$-modules $I_t / (1 - Ad (\Phi)^m)I_t \stackrel{\sim}{\rightarrow} \mathfrak{f}_m^*$ where $Ad$ denotes the adjoint action.  Since $\hat{\sigma}$ has order $n$, we have $\mathrm{Hom}_{\Phi, \hat{\sigma}}(\mathfrak{f}_n^*, \hat{T}) \cong \mathrm{Hom}_{Ad(\Phi), \hat{\sigma}}(I_t, \hat{T})$.  Therefore, the map $s \mapsto \chi_s$ is a canonical bijection $$\mathrm{Hom}_{Ad(\Phi), \hat{\sigma}}(I_t, \hat{T}) \stackrel{\sim}{\rightarrow} \mathrm{Hom}(\mathbb{T}^{\Phi_{\sigma}}, \mathbb{C}^*)$$
Moreover, we have an isomorphism $${}^0 T^{\Phi_{\sigma}} \times X^{\sigma} \stackrel{\sim}{\rightarrow} T^{\Phi_{\sigma}}$$ $$(\gamma, \lambda) \mapsto \gamma \lambda(\varpi)$$ where ${}^0 T$ is the group of $\mathfrak{o}_{F^u}$-points of $\mathbf{T}$.

Finally, note that $\hat{T} / (1 - \hat{\sigma}) \hat{T}$ is the character group of $X^{\sigma}$, whereby $\tau \in \hat{T} / (1 - \hat{\sigma}) \hat{T}$ corresponds to $\chi_{\tau} \in \mathrm{Hom}(X^{\sigma}, \mathbb{C}^*),$ where $ \chi_{\tau}(\lambda) := \lambda(\tau)$.  Therefore, we have a canonical bijection between $\hat{T}$-conjugacy classes of admissible homomorphisms $\phi : W_t \rightarrow {}^L T_{\sigma}$ and depth-zero characters

\begin{equation}
\chi_{\phi} := \chi_s \otimes \chi_{\tau} \in \mathrm{Irr}(T^{\Phi_{\sigma}}) \ \ \label{chitau}
\end{equation}

\noindent where $s := \phi|_{I_t}$, $\phi(\Phi) = \hat{\sigma} \ltimes \tau$, and where we have inflated $\chi_s$ to ${}^0 T^{\Phi_{\sigma}}$.

\

To get the depth-zero $L$-packet associated to $\phi$, one implements the component group $$\mathrm{Irr}(C_{\phi}) \cong [X / (1 - w \theta) X]_{\mathrm{tor}}$$ as follows. We set $X_w$ to be the preimage of $[X / (1 - w \theta) X]_{\mathrm{tor}}$ in $X$.  To $\lambda \in X_w$, DeBacker and Reeder associate a 1-cocycle $u_{\lambda}$, hence a twisted Frobenius $\Phi_{\lambda} = Ad(u_{\lambda}) \circ \Phi$.  Moreover, to $\lambda$, they associate a facet $J_{\lambda}$, and hence a parahoric subgroup $G_{\lambda}$ associated to $J_{\lambda}$.  Let $\mathbb{G}_{\lambda} := G_{\lambda} / G_{\lambda}^+$.  Let $W_{\lambda}$ be the subgroup of $W^o$ generated by reflections in the hyperplanes containing $J_{\lambda}$.  Then to $\lambda$, DeBacker--Reeder also associate an element $w_{\lambda} \in W_{\lambda}$. Fix once and for all a lift $\dot{w}$ of $w$ to $N$.  Using this lift, DeBacker and Reeder also associate a lift $\dot{w}_{\lambda}$ of $w_{\lambda}$ to $N$. By Lang's theorem, there exists $p_{\lambda} \in G_{\lambda}$ such that $p_{\lambda}^{-1} \Phi_{\lambda} (p_{\lambda}) = \dot{w}_{\lambda}$.  We then define $T_{\lambda} := Ad(p_{\lambda}) T$, and set $\chi_{\lambda} := \chi_{\phi} \circ Ad(p_{\lambda})^{-1}$.  Since $\chi_{\lambda}$ is depth-zero, its restriction to ${}^0 T_{\lambda}^{\Phi_{\lambda}}$ factors through a character $\chi_{\lambda}^0$ of $\mathbb{T}_{\lambda}^{\Phi_{\lambda}}$, where $\mathbb{T}_{\lambda}^{\Phi_{\lambda}}$ is the projection of ${}^0 T^{\Phi_{\lambda}}$ in $\mathbb{G}_{\lambda}$.  Therefore, $\chi_{\lambda}^0$ gives rise to an irreducible cuspidal Deligne-Lusztig representation $\kappa_{\lambda}^0$ of $\mathbb{G}_{\lambda}^{\Phi_{\lambda}}$.  Inflate $\kappa_{\lambda}^0$ to a representation of $G_{\lambda}^{\Phi_{\lambda}}$, and define an extension to $Z^{\Phi_{\lambda}} G_{\lambda}^{\Phi_{\lambda}}$ by $$\kappa_{\lambda} := \chi_{\lambda} \otimes \kappa_{\lambda}^0$$ where $Z$ denotes the center of $G$.  This makes sense since $(Z \cap G_{\lambda})^{\Phi_{\lambda}}$ acts on $\kappa_{\lambda}^0$ via the restriction of $\chi_{\lambda}^0$.  Finally, form the representation $$\pi_{\lambda} := \mathrm{Ind}_{Z^{\Phi_{\lambda}} G_{\lambda}^{\Phi_{\lambda}}}^{G^{\Phi_{\lambda}}} \kappa_{\lambda}$$ where Ind denotes smooth induction.  Then DeBacker--Reeder construct a packet $\Pi(\phi)$ of representations on the pure inner forms of $G$, parameterized by $\mathrm{Irr}(C_{\phi})$, using the above construction, where $C_{\phi}$ is the component group of $\phi$.

\section{Existing Description of the Tame Local Langlands Correspondence for $GL(\ell,F)$, $\ell$ a prime}\label{existingdescription}

In this chapter, we describe the construction of the tame local Langlands correspondence for $GL(\ell,F)$ as explained in \cite{moy}, where $\ell$ is a prime.

\subsection{Depth-zero supercuspidal representations of $GL(\ell,F)$, $\ell$ a prime}\label{depthzeroclassicalconstruction}

Let $(E/F, \chi)$ be an admissible pair where $\chi$ has level 0 and $E/F$ has degree $\ell$.  By definition of admissible pair, this implies that $E/F$ is unramified, and the residue field of $E$ is $\mathfrak{f}_{\ell}$.  Since $\chi|_{1 + \mathfrak{p}_E} = 1$, $\chi|_{\mathfrak{o}_E^*}$ is the inflation of the character $\chi_o$ of $\mathfrak{f}_{\ell}^*$.  By the theory of finite groups of Lie type, this character gives rise to an irreducible cuspidal representation $\lambda'$ of $GL(\ell,\mathfrak{f})$, which is the irreducible cuspidal Deligne-Lusztig representation corresponding to the elliptic torus $\mathfrak{f}_{\ell}^* \subset GL(\ell,\mathfrak{f})$ and the character $\chi_o$ of $\mathfrak{f}_{\ell}^*$.  Let $\lambda$ be the inflation of $\lambda'$ to $GL(\ell, \mathfrak{o}_F)$.  We may extend $\lambda$ to a representation $\Lambda$ of $K(F) := F^* GL(\ell, \mathfrak{o}_F)$ by setting $\Lambda|_{F^*} = \chi|_{F^*}$, and then induce the resulting representation to $G(F) = GL(\ell,F)$.  Set $$\pi_{\chi} = cInd_{K(F)}^{G(F)} \Lambda$$ where $cInd$ denotes compact induction.  Let $\mathbb{P}_{\ell}(F)_0$ be the subset of admissible pairs $(E/F, \chi)$ such that $\chi$ has level zero and $\mathbb{A}_{\ell}^0(F)_0$ be the subset of depth-zero supercuspidal representations of $GL(\ell,F)$.

\begin{proposition}\label{depthzerogl2}
Suppose that $p \neq \ell$.  The map (E/F, $\chi$) $\mapsto \pi_{\chi}$ induces a bijection $$\mathbb{P}_{\ell}(F)_0 \rightarrow \mathbb{A}_{\ell}^0(F)_0$$
\end{proposition}

\proof

See \cite{moy}.

\qed

\subsection{Positive depth supercuspidal representations of $GL(\ell,F)$, $\ell$ a prime}

In this section we recall the parameterization of the positive depth supercuspidal representions via admissible pairs, following \cite{moy}.  Let $\mathbb{A}_{\ell}^0(F)^+$ denote the set of all positive depth irreducible supercuspidal representations of $GL(\ell,F)$, and let $\mathbb{P}_{\ell}(F)^+$ denote the set of all admissible pairs $(E/F, \chi) \in \mathbb{P}_{\ell}(F)$ such that $\chi$ has positive level.

\begin{proposition}\label{bijection2}
Suppose that $p \neq \ell$.  There is a map $(E/F, \chi)$ $\mapsto \pi_{\chi}$ that induces a bijection $$\mathbb{P}_{\ell}(F)^+ \rightarrow \mathbb{A}_{\ell}^0(F)^+$$
\end{proposition}

\proof
See \cite{moy}.
\qed

\subsection{Langlands parameters}\label{weilparameters}

Let $\mathbb{G}_{\ell}^0(F)$ be the set of equivalence classes of irreducible smooth $\ell$-dimensional representations of $W_F$. Recall that there is a local Artin reciprocity isomorphism given by $W_E^{ab} \cong E^*$.  Then, if $(E/F, \chi) \in \mathbb{P}_{\ell}(F)$, $\chi$ gives rise to a character of $W_E^{ab}$, which we can pullback to a character, also denoted $\chi$, of $W_E$.  We can then form the induced representation $\phi(\chi) := Ind_{W_E}^{W_F} \chi$ of $W_F$.

\begin{theorem}\label{G_2(F)}
Suppose $p \neq \ell$.  If $(E/F, \chi) \in \mathbb{P}_{\ell}(F)$, the representation $\phi(\chi)$ of $W_F$ is irreducible.  The map $(E/F, \chi) \mapsto \phi(\chi)$ induces a bijection $$\mathbb{P}_{\ell}(F) \rightarrow \mathbb{G}_{\ell}^0(F)$$
\end{theorem}

\proof
See \cite{moy}.
\qed

For the next theorem, we will need to associate to any admissible pair $(E/F, \chi) \in \mathbb{P}_{\ell}(F)$ a specific character $\Delta_{\chi}$ of $E^*$.  We will not define $\Delta_{\chi}$ in general, but only for the cases that we need in this paper.  For the general definition of $\Delta_{\chi}$ associated to any admissible pair $(E/F, \chi) \in \mathbb{P}_{\ell}(F)$, see \cite{moy}.

\begin{definition}
If $(E/F, \chi)$ is an admissible pair in which $E/F$ is quadratic and unramified, define $\Delta_{\chi}$ to be the unique quadratic unramified character of $E^*$.  If $(E/F, \chi)$ is an admissible pair in which $E/F$ is degree $\ell$ and unramified, where $\ell$ is an odd prime, then define $\Delta_{\chi}$ to be the trivial character of $E^*$.
\end{definition}

\begin{theorem}{$\mathbf{Tame \ Local \ Langlands \ Correspondence}$}{\label{tamellc}}

Suppose $p \neq \ell$.  For $\phi \in \mathbb{G}_{\ell}^0(F)$, define $\pi(\phi) = \pi_{\chi \Delta_{\chi}}$ in the notation of Propositions \ref{depthzerogl2} and \ref{bijection2} for any $(E/F, \chi) \in \mathbb{P}_{\ell}(F)$ such that $\phi \cong \phi(\chi)$.  The map $$\pi : \mathbb{G}_{\ell}^0(F) \rightarrow \mathbb{A}_{\ell}^0(F)$$ is the local Langlands correspondence for supercuspidal representations of $GL(\ell,F)$.
\end{theorem}

\proof
See \cite{moy}.
\qed

\section{The case of $GL(\ell,F)$, $\ell$ a prime}\label{glellf}

For section \ref{glellf} and section \ref{glellfcharactertorus}, we consider the group $\mathbf{G}(F) = GL(\ell,F)$, where $\ell$ is prime.  We will show that the conjectural correspondence of \cite{debackerreeder} agrees with the local Langlands correspondence for $GL(\ell,F)$ given in Section \ref{existingdescription}.

Let $\phi : W_F \rightarrow {}^L G$ be a TRSELP for $\mathbf{G}(F) = GL(\ell,F)$.  This is equivalent to an irreducible admissible $\phi : W_F \rightarrow GL(\ell,\mathbb{C})$ that is trivial on the wild inertia group.  By Section \ref{weilparameters}, have that $\phi = Ind_{W_E}^{W_F}(\chi)$ for some admissible pair $(E/F, \chi)$, where $\chi$ has level zero and $E/F$ is degree $\ell$ and unramified.  We will need the relative Weil group $W_{E/F} := W_F / [W_E, W_E]^c$ (see \cite[Chapter 1]{tate}), where $c$ denotes closure and $[W_E, W_E]$ denotes the commutator subgroup of $W_E$.  The representation $\phi = Ind_{W_E}^{W_F} (\chi)$ factors through $W_{E/F}$, since $\phi|_{W_E} = \chi \oplus \chi^{\xi} \oplus ... \oplus \chi^{\xi^{\ell-1}}$.

We begin by calculating the character $\chi_{\phi}$ from (\ref{chitau}).  Note that ${}^L G = \langle \hat{\theta} \rangle \times GL(\ell,\mathbb{C})$.  $\hat{T}$ is the diagonal maximal torus in $\hat{G} = GL(\ell,\mathbb{C})$, and after conjugation, we may assume $\phi(I_F) \subset \hat{T}$.  Moreover, $\phi(\Phi) = \hat{\theta} f$ for some $f \in \hat{N}$ such that $\hat{w}$ is a cycle of length $\ell$ in the Weyl group $S_{\ell}$, the symmetric group on $\ell$ letters.  The reason for this requirement on the Weyl group element is that $\phi$ is TRSELP and so in particular it is elliptic.  In particular, ellipticity is equivalent to requiring that the image of $\phi$ is not contained in any proper Levi subgroup of ${}^L G$ (see \cite[section 3.4]{debackerreeder}).  After conjugating the TRSELP by a permutation matrix in $N_{\hat{G}}(\hat{T})$, we may assume without loss of generality that $\hat{w} = (1 \ 2 \ 3 \ \cdots \ \ell) \in S_{\ell}$ since all cycles of length $\ell$ are conjugate in $S_{\ell}$.  Note that this choice implies that $w = (1 \ 2 \ 3 \ \cdots \ \ell) \in S_{\ell}$.  The arguments in the remainder of the paper are the same for all other allowable choices of $\hat{w}$.

Let us first calculate $\chi_s$, where $s := \phi|_{I_t}$ (recall again that $\phi|_{I_F^+} \equiv 1$, so $\phi|_{I_F}$ factors to $I_t$).

\begin{proposition}\label{whatevers2}
Let $\phi = Ind_{W_E}^{W_F}(\chi)$ and set $s = \phi|_{I_t}$, where $(E/F, \chi)$ is an admissible pair as above. Then the isomorphism $$\mathrm{Hom}_{Ad(\Phi), \hat{\sigma}}(I_t, \hat{T}) \stackrel{\sim}{\rightarrow} \mathrm{Hom}_{\Phi, \hat{\sigma}}(\mathfrak{f}_{\ell}^*, \hat{T})$$ sends $s$ to $\tilde\beta_s$, where
\[\tilde\beta_s(x) = \left( \begin{array}{ccccc}
\chi_o(x) & 0 & 0 & \cdots & 0\\
0 & \chi_o^{\xi}(x) & 0 & \cdots & 0\\
0 & 0 & \chi_o^{\xi^2}(x) & \cdots & 0\\
\vdots & \vdots & \vdots & \ddots & \vdots\\
0 & 0 & 0 & \cdots & \chi_o^{\xi^{\ell-1}}(x)
\end{array} \right) \]
\end{proposition}

\proof
Since $\hat{\sigma}$ has order $\ell$, $s \in \mathrm{Hom}_{Ad(\Phi), \hat{\sigma}}(I_t, \hat{T})$ is trivial on $(1 - Ad(\Phi)^{\ell}) I_t$, so factors to $I_t / (1 - Ad(\Phi)^{\ell}) I_t$.  We first note that the isomorphisms

\begin{equation*}
I_t \cong \lim_{\stackrel{\leftarrow}{m}}
\mathfrak{f}_m^*, \ \ \ \ \ I_t / (1 - Ad(\Phi)^{\ell}) I_t \stackrel{\sim}{\rightarrow} \mathfrak{f}_{\ell}^*
\end{equation*}

\noindent are induced by local Artin Reciprocity (see \cite[Chapter 5]{reeder}).  Moreover, the map $\mathrm{Hom}_{Ad(\Phi), \hat{\sigma}}(I_t, \hat{T}) \rightarrow \mathrm{Hom}_{\Phi, \hat{\sigma}}(\mathfrak{f}_{\ell}^*, \hat{T})$ comes from the following diagram:

\begin{diagram}
I_t  & \rTo^s & \hat{T} \\
\dTo & \ruTo  &  \\
I_t / (1 - Ad(\Phi)^{\ell}) I_t \\
\uTo^{\wr} \\
\mathfrak{f}_{\ell}^*
\end{diagram}

\noindent Recall that $\phi$ factors through $W_{E/F}$.  Therefore, we also have the following commutative diagram.

\begin{diagram}
W_F  & \rTo^{\phi} & GL(\ell,\mathbb{C}) \\
\dTo & \ruTo^{\beta}  &  \\
W_{E/F}
\end{diagram}

\noindent It is a fact that $W_{E/F}$ is an extension of $Gal(E/F)$ by $E^*$, and can be described by generators and relations as follows.  The generators are $\{z \in E^* \}$ and an element $j$ where $j \in W_{E/F}$ satisfies $j^{\ell} = \varpi$ and $jzj^{-1} = \xi(z)$.  Then the map $W_F \rightarrow W_{E/F}$ sends $I_F$ to $\mathfrak{o}_E^*$ and $\Phi$ to $j$.

Let us calculate the map $\beta$.  Consider the canonical sequence $$1 \rightarrow W_E / [W_E, W_E]^c \rightarrow W_F / [W_E, W_E]^c \rightarrow W_F / W_E \cong Gal(E/F) \rightarrow 1$$
Recall that $\phi$ is trivial on $[W_E, W_E]^c$.  To calculate $\beta|_{E^*}$, it suffices to calculate $\phi|_{W_E}$ since $W_E / [W_E, W_E]^c \cong E^*$ by Artin reciprocity.  But $\phi|_{W_E} = \chi \oplus \chi^{\xi} \oplus ... \oplus \chi^{\xi^{\ell-1}}$.  Therefore,

\[\beta(t) = \left( \begin{array}{ccccc}
\chi(t) & 0 & 0 & \cdots & 0\\
0 & \chi^{\xi}(t) & 0 & \cdots & 0\\
0 & 0 & \chi^{\xi^2}(t) & \cdots & 0\\
\vdots & \vdots & \vdots & \ddots & \vdots\\
0 & 0 & 0 & \cdots & \chi^{\xi^{\ell-1}}(t)
\end{array} \right) \]

Moreover, since $\phi$ is irreducible, we have that $\beta(j) \in N_{GL(\ell,\mathbb{C})}(\hat{T})$ represents $\hat{w}$.

Since $\phi|_{I_F^+} \equiv 1$, we have that $\beta|_{1 + \mathfrak{p}_E} \equiv 1$, so $\beta|_{\mathfrak{o}_E^*}$ factors to a map $\tilde\beta_s : \mathfrak{f}_{\ell}^* \rightarrow GL(\ell,\mathbb{C})$, given by

\[\tilde\beta_s(x) = \left( \begin{array}{ccccc}
\chi_o(x) & 0 & 0 & \cdots & 0\\
0 & \chi_o^{\xi}(x) & 0 & \cdots & 0\\
0 & 0 & \chi_o^{\xi^2}(x) & \cdots & 0\\
\vdots & \vdots & \vdots & \ddots & \vdots\\
0 & 0 & 0 & \cdots & \chi_o^{\xi^{\ell-1}}(x)
\end{array} \right) \ \ \forall x \in \mathfrak{f}_{\ell}^* \]
\qed

\begin{proposition}\label{bigcomposite4}
Let $\phi = Ind_{W_E}^{W_F}(\chi)$ and set $s = \phi|_{I_t}$ as above.  Then the composite isomorphism
$$\mathrm{Hom}_{Ad(\Phi), \hat{\sigma}}(I_t, \hat{T}) \stackrel{\sim}{\rightarrow} \mathrm{Hom}_{\Phi, \hat{\sigma}}(\mathfrak{f}_{\ell}^*, \hat{T}) \stackrel{\sim}{\rightarrow} \mathrm{Hom}_{\Phi_{\sigma}, Id}(X \otimes \mathfrak{f}_{\ell}^*, \mathbb{C}^*)$$ $$ \stackrel{\sim}{\rightarrow} \mathrm{Hom}_{\Phi_{\sigma}, Id}(\mathbb{T}^{\Phi_{\sigma}^{\ell}}, \mathbb{C}^*)$$
sends $s$ to ${}^{\ell} \chi_o^E$, where ${}^{\ell} \chi_o^E(x_1,x_2, \cdots, x_{\ell}) := \chi_o(x_1) \chi_o^{\xi}(x_2) \cdots \chi_o^{\xi^{\ell-1}}(x_{\ell})$ and where $s = \phi|_{I_t}, \phi = Ind_{W_E}^{W_F}(\chi)$.
\end{proposition}

\proof
The composite isomorphism $$\mathrm{Hom}_{\Phi, \hat{\sigma}}(\mathfrak{f}_{\ell}^*, \hat{T}) \stackrel{\sim}{\rightarrow} \mathrm{Hom}_{\Phi_{\sigma}, Id}(X \otimes \mathfrak{f}_{\ell}^*, \mathbb{C}^*) \stackrel{\sim}{\rightarrow} \mathrm{Hom}_{\Phi_{\sigma}, Id}(\mathbb{T}^{\Phi_{\sigma}^{\ell}}, \mathbb{C}^*)$$ is given by $$\tilde\alpha \mapsto \{ \lambda(x) \mapsto \lambda(\tilde\alpha(x)) \}$$
where $x \in \mathfrak{f}_{\ell}^*$ and $\lambda \in X = X_*(T) = X^*(\hat{T})$.  $\mathbb{T}$ splits over $\mathfrak{f}_{\ell}$ and note that $\mathbb{T}^{\Phi_{\sigma}^{\ell}} \cong \mathfrak{f}_{\ell}^* \times ... \times \mathfrak{f}_{\ell}^*$. Then, it is easy to see that under this composite isomorphism, $\tilde\beta_s$ (where $\tilde\beta_s$ is as in Proposition \ref{whatevers2}) maps to the homomorphism

\[ (x_1, x_2, ..., x_{\ell})  \mapsto \chi_o(x_1) \chi_o^{\xi}(x_2) \cdots \chi_o^{\xi^{\ell-1}}(x_{\ell}) \ \forall x_1, x_2, ..., x_{\ell} \in \mathfrak{f}_{\ell}^*. \]

\noindent by considering the standard basis of cocharacters of $X$.
\qed

\begin{proposition}\label{whatevers3}
The isomorphism $\mathrm{Hom}_{\Phi_{\sigma}, Id}(\mathbb{T}^{\Phi_{\sigma}^{\ell}}, \mathbb{C}^*) \stackrel{\sim}{\rightarrow} \mathrm{Hom}(\mathbb{T}^{\Phi_{\sigma}}, \mathbb{C}^*)$ is given by

\begin{equation}
\mathrm{Hom}_{\Phi_{\sigma}, Id}(\mathbb{T}^{\Phi_{\sigma}^{\ell}}, \mathbb{C}^*) \stackrel{\sim}{\rightarrow} \mathrm{Hom}(\mathbb{T}^{\Phi_{\sigma}}, \mathbb{C}^*)  \ \ \label{chichi'ell}
\end{equation}

$$\Lambda \mapsto \Lambda'$$

\noindent where $\Lambda'(a) := \Lambda((x_1, x_2, \cdots, x_{\ell}))$ for any $(x_1, x_2, \cdots, x_{\ell}) \in \mathfrak{f}_{\ell}^* \times \mathfrak{f}_{\ell}^* \times \cdots \times \mathfrak{f}_{\ell}^*$ such that \\ $x_1 x_2^{q^{\ell-1}} x_3^{q^{\ell-2}} \cdots x_{\ell}^q = a$, where $ a \in \mathfrak{f}_{\ell}^*$.
\end{proposition}

\proof
Recall that the isomorphism $\mathrm{Hom}_{\Phi_{\sigma}, Id}(\mathbb{T}^{\Phi_{\sigma}^{\ell}}, \mathbb{C}^*) \stackrel{\sim}{\rightarrow} \mathrm{Hom}(\mathbb{T}^{\Phi_{\sigma}}, \mathbb{C}^*)$ is abstractly given by

\begin{equation*}
\mathrm{Hom}_{\Phi_{\sigma}, Id}(\mathbb{T}^{\Phi_{\sigma}^{\ell}}, \mathbb{C}^*) \stackrel{\sim}{\rightarrow} \mathrm{Hom}(\mathbb{T}^{\Phi_{\sigma}}, \mathbb{C}^*)
\end{equation*}

$$\Lambda \mapsto \Lambda'$$

\noindent where $\Lambda'(a) := \Lambda((x_1, x_2, ..., x_{\ell}))$ for any $(x_1, x_2, ..., x_{\ell}) \in \mathfrak{f}_{\ell}^* \times \mathfrak{f}_{\ell}^* \times \cdots \times \mathfrak{f}_{\ell}^*$ such that \\ $N_{\sigma}((x_1, x_2, ..., x_{\ell})) = a$.

We need some preliminaries.  First note that $\Phi_{\sigma}((x_1, x_2, ..., x_{\ell})) = w \Phi^{-1}((x_1, x_2, ..., x_{\ell})) = w (x_1^q, x_2^q, ..., x_{\ell}^q) = (x_{\ell}^q, x_1^q, x_2^q, ..., x_{\ell-1}^q)$.  If we make the identification of $\mathbb{T}^{\Phi_{\sigma}^{\ell}}$ with tuples $(x_1, x_2, ..., x_{\ell}) \in \mathfrak{f}_{\ell}^* \times \mathfrak{f}_{\ell}^* \times \cdots \times \mathfrak{f}_{\ell}^*$, then we have that since we made our choice of $w = (1 \ 2 \ 3 \ \cdots \ \ell) \in S_{\ell}$, we get $\mathbb{T}^{\Phi_{\sigma}} = \{ (x_1, x_2, ..., x_{\ell}) \in \mathfrak{f}_{\ell}^* \times \mathfrak{f}_{\ell}^* \times \cdots \times \mathfrak{f}_{\ell}^* : (x_{\ell}^q, x_1^q, x_2^q, \cdots, x_{\ell-1}^q) = (x_1, x_2, ..., x_{\ell}) \} = \{(x_1, x_1^q, x_1^{q^2}, \cdots, x_1^{q^{\ell-1}}) : x_1 \in \mathfrak{f}_{\ell} \}$.  If $(x_1, x_2, ..., x_{\ell}) \in \mathfrak{f}_{\ell}^* \times \mathfrak{f}_{\ell}^* \times \cdots \times \mathfrak{f}_{\ell}^* = \mathbb{T}^{\Phi_{\sigma}^{\ell}}$, then $$N_{\sigma}((x_1, x_2, ..., x_{\ell})) = (x_1, x_2, ..., x_{\ell}) \Phi_{\sigma}((x_1, x_2, ..., x_{\ell})) \Phi_{\sigma}^2((x_1, x_2, ..., x_{\ell})) \cdots \Phi_{\sigma}^{\ell-1}((x_1, x_2, ..., x_{\ell})) =$$
$$(x_1, x_2, ..., x_{\ell}) (x_{\ell}^q, x_1^q, x_2^q, \cdots, x_{\ell-1}^q) (x_{\ell-1}^{q^2}, x_{\ell}^{q^2}, x_1^{q^2}, \cdots, x_{\ell-2}^{q^2}) \cdots (x_2^{q^{\ell-1}}, x_3^{q^{\ell-1}}, \cdots, x_{\ell}^{q^{\ell-1}}, x_1^{q^{\ell-1}}) = $$ $$(x_1 x_2^{q^{\ell-1}} x_3^{q^{\ell-2}} \cdots x_{\ell}^q, x_2 x_3^{q^{\ell-1}} x_4^{q^{\ell-2}} \cdots x_1^q, \cdots, x_{\ell} x_1^{q^{\ell-1}} x_2^{q^{\ell-2}} \cdots x_{\ell-1}^q)$$
Therefore, $N_{\sigma} : \mathbb{T}^{\Phi_{\sigma}^{\ell}} \rightarrow \mathbb{T}^{\Phi_{\sigma}}$ is the map $$ (x_1, x_2, ..., x_{\ell}) \mapsto x_1 x_2^{q^{\ell-1}} x_3^{q^{\ell-2}} \cdots x_{\ell}^q, \ \ \forall (x_1, x_2, ..., x_{\ell}) \in \mathfrak{f}_{\ell}^* \times \mathfrak{f}_{\ell}^* \times \cdots \times \mathfrak{f}_{\ell}^*$$
\qed

We now need to obtain a character of ${}^0 T^{\Phi_{\sigma}}$ from a character of $\mathbb{T}^{\Phi_{\sigma}}$.  In our case, ${}^0 T^{\Phi_{\sigma}} = \mathfrak{o}_{E}^*$, which has a canonical projection map $\mathfrak{o}_{E}^* = {}^0 T^{\Phi_{\sigma}} \xrightarrow{\eta} \mathbb{T}^{\Phi_{\sigma}} = \mathfrak{f}_{\ell}^*$.  Then, given $\zeta \in \mathrm{Hom}(\mathbb{T}^{\Phi_{\sigma}}, \mathbb{C}^*)$, we obtain a character $\mu$ of ${}^0 T^{\Phi_{\sigma}} = \mathfrak{o}_{E}^*$ given by $\mu(z) := \zeta(\eta(z)), z \in \mathfrak{o}_{E}^*$.  Let us more explicitly calculate such a $\mu$, given some $\Lambda' \in \mathrm{Hom}(\mathbb{T}^{\Phi_{\sigma}}, \mathbb{C}^*)$ that comes from  $\Lambda \in \mathrm{Hom}_{\Phi_{\sigma}, Id}(\mathbb{T}^{\Phi_{\sigma}^{\ell}}, \mathbb{C}^*)$ as in (\ref{chichi'ell}).  Let $z \in \mathfrak{o}_{E}^*$.  Then $\mu(z) = \Lambda'(\eta(z)) = \Lambda((\eta(z),1, 1, \cdots, 1))$, by Proposition \ref{whatevers3}.

We may now calculate the character $\chi_{s}$ that arises from $\phi$, where $s = \phi|_{I_t}$ and $\phi = Ind_{W_E}^{W_F}(\chi)$.  The above analysis and Proposition \ref{bigcomposite4} shows that $\chi_{s}(z) = {}^{\ell} \chi_o^E((\eta(z),1, 1, \cdots, 1)) = \chi_o(\eta(z)) = \chi(z)$, where $z \in \mathfrak{o}_{E}^*$.  It remains to compute $\chi_{\tau}$.  First note that if we make the identification $X = \mathbb{Z} \times \mathbb{Z} \times \cdots \times \mathbb{Z}$, then $X^{\sigma} = \{(k,k, \cdots, k) : k \in \mathbb{Z} \}$.  Let $\lambda_{(k,k, \cdots, k)} \in X^{\sigma}$ denote the character of $\hat{T}$ corresponding to $(k,k, \cdots, k) \in \mathbb{Z} \times \mathbb{Z} \times \cdots \times \mathbb{Z}$ via this identification.

\begin{proposition}\label{chitaugl2}
Let $\ell = 2$.  The character $\chi_{\tau}$ is given by $$\chi_{\tau} : X^{\sigma} \rightarrow \mathbb{C}^*$$ $$\lambda_{(k,k)} \mapsto (-\chi(\varpi))^k$$
\end{proposition}

\proof
Note that $\hat{\theta} = 1$ and $\hat{G}' = SL(2,\mathbb{C})$, so $\tau$ is any element whose class in $\hat{T} / (1 - \hat{\sigma}) \hat{T}$ corresponds to the image of $f$ in $GL(2,\mathbb{C}) / SL(2,\mathbb{C})$ under the bijection $$\hat{T} / (1 - \hat{\sigma}) \hat{T} \stackrel{\sim}{\rightarrow} GL(2,\mathbb{C}) / SL(2,\mathbb{C})$$ as in (\ref{bijectionfortau}).  We thus need to compute $f$ first.

Recall that $\phi(\Phi) = \beta(j)$, where $\beta$ is as in the proof of Proposition \ref{whatevers2}.  Recall that since $\phi$ is irreducible, then $\beta(j) = \mat{0}{a}{b}{0}$ for some $a,b \in \mathbb{C}^*$.  After conjugation by $\hat{G}$, we may assume that $b = 1$.  But since $j^2 = \varpi$, we have $$\mat{\chi(\varpi)}{0}{0}{\chi^{\xi}(\varpi)} = \beta(\varpi) = \beta(j^2) = \beta(j)^2 = \mat{a}{0}{0}{a}$$  Therefore, $a = \chi(\varpi)$ and so $\beta(j) = \mat{0}{\chi(\varpi)}{1}{0}$, and we may take $f = \mat{0}{\chi(\varpi)}{1}{0}$.

We now note that the bijection $$\hat{T} / (1 - \hat{\sigma}) \hat{T} \stackrel{\sim}{\rightarrow} \hat{G}_{ab} / (1 - \hat{\theta}) \hat{G}_{ab}$$ is induced by the inclusion $\hat{T} \hookrightarrow \hat{G}$ (see \cite[section 4.3]{debackerreeder}).  Now, we have that $$\left[\mat{-\chi(\varpi)}{0}{0}{1} \right] = \left[\mat{0}{\chi(\varpi)}{1}{0} \right]  \in GL(2,\mathbb{C}) / SL(2,\mathbb{C})$$ since $$\mat{0}{\chi(\varpi)}{1}{0}^{-1} \mat{-\chi(\varpi)}{0}{0}{1} = \mat{0}{1}{-1}{0} \in SL(2,\mathbb{C})$$
Therefore, since $\mat{-\chi(\varpi)}{0}{0}{1} \in \hat{T}$, we may set $\tau = \mat{-\chi(\varpi)}{0}{0}{1}$.  Then $\chi_{\tau} : X^{\sigma} \rightarrow \mathbb{C}^*$ is given by $\chi_{\tau}(\lambda_{(k,k)}) = \lambda_{(k,k)}(\tau) = \lambda_{(k,k)}\mat{-\chi(\varpi)}{0}{0}{1} = (-\chi(\varpi))^k$.
\qed

\begin{proposition}\label{chitauglelldepthzero}
Let $\ell$ be an odd prime.  The character $\chi_{\tau}$ is given by $$\chi_{\tau} : X^{\sigma} \rightarrow \mathbb{C}^*$$ $$\lambda_{(k,k, \cdots, k)} \mapsto \chi(\varpi)^k$$
\end{proposition}

\proof
Note that $\hat{\theta} = 1$ and $\hat{G}' = SL(\ell,\mathbb{C})$, so $\tau$ is any element whose class in $\hat{T} / (1 - \hat{\sigma}) \hat{T}$ corresponds to the image of $f$ in $GL(\ell,\mathbb{C}) / SL(\ell,\mathbb{C})$ under the bijection $$\hat{T} / (1 - \hat{\sigma}) \hat{T} \stackrel{\sim}{\rightarrow} GL(\ell,\mathbb{C}) / SL(\ell,\mathbb{C})$$ as in (\ref{bijectionfortau}).  We thus need to compute $f$ first.

Recall that $\phi$ factors through $W_{E/F}$, and we have the commutative diagram

\begin{diagram}
W_F  & \rTo^{\phi} & GL(\ell,\mathbb{C}) \\
\dTo & \ruTo^{\beta}  &  \\
W_{E/F}
\end{diagram}

\noindent From Proposition \ref{whatevers2}, we have $\phi(\Phi) = \beta(j)$.  To compute $\beta(j)$, recall that because of our choice of $\hat{w}$, we have
\[ \beta(j) = \left( \begin{array}{ccccc}
0 & 0 & 0 & \cdots & a_1\\
a_2 & 0 & 0 & \cdots & 0\\
0 & a_3 & 0 & \cdots & 0\\
\vdots & \vdots & \vdots & \vdots & \vdots\\
0 & 0 & \cdots & a_{\ell} & 0
\end{array} \right) \]

\noindent for some $a_i \in \mathbb{C}^*$.  After conjugation the Langlands parameter by an element in $\hat{G}$ of the form

\[ \left( \begin{array}{cccccc}
0 & x_2 & 0 & 0 & \cdots & 0\\
0 & 0 & x_3 & 0 & \cdots & 0\\
0 & 0 & 0 & x_4 & \cdots & \vdots\\
\vdots & \vdots & \vdots & \vdots & \ddots & 0 \\
0 & 0 & 0 & 0 & 0 & x_{\ell}\\
x_1 & 0 & 0 & \cdots & 0 & 0
\end{array} \right), \]

we may assume that $a_2 = a_3 = \cdots = a_{\ell} = 1$.   Therefore,

\[\left( \begin{array}{ccccc}
\chi(\varpi) & 0 & 0 & \cdots & 0\\
0 & \chi^{\xi}(\varpi) & 0 & \cdots & 0\\
0 & 0 & \chi^{\xi^2}(\varpi) & \cdots & 0\\
\vdots & \vdots & \vdots & \vdots & \vdots\\
0 & 0 & 0 & \cdots & \chi^{\xi^{\ell-1}}(\varpi)
\end{array} \right) = \beta(\varpi) = \beta(j^{\ell}) = \]

\[\beta(j)^{\ell} =
\left( \begin{array}{ccccc}
a_1 & 0 & 0 & \cdots & 0\\
0 & a_1 & 0 & \cdots & 0\\
0 & 0 & a_1 & \cdots & 0\\
\vdots & \vdots & \vdots & \vdots & \vdots\\
0 & 0 & 0 & \cdots & a_1
\end{array} \right)\]

\noindent Therefore, $a_1 = \chi(\varpi)$ and so we may take
\[ f = \left( \begin{array}{ccccc}
0 & 0 & 0 & \cdots & \chi(\varpi)\\
1 & 0 & 0 & \cdots & 0\\
0 & 1 & 0 & \cdots & 0\\
\vdots & \vdots & \vdots & \vdots & \vdots\\
0 & 0 & \cdots & 1 & 0
\end{array} \right) \]

Now, we have that

\[ \left[ \left( \begin{array}{ccccc}
\chi(\varpi) & 0 & 0 & \cdots & 0\\
0 & 1 & 0 & \cdots & 0\\
0 & 0 & 1 & \cdots & 0\\
\vdots & \vdots & \vdots & \vdots & \vdots\\
0 & 0 & \cdots & 0 & 1
\end{array} \right) \right]
= \left[ \left( \begin{array}{ccccc}
0 & 0 & 0 & \cdots & \chi(\varpi)\\
1 & 0 & 0 & \cdots & 0\\
0 & 1 & 0 & \cdots & 0\\
\vdots & \vdots & \vdots & \vdots & \vdots\\
0 & 0 & \cdots & 1 & 0
\end{array} \right) \right] \in GL(\ell,\mathbb{C}) / SL(\ell,\mathbb{C})\]

since

\[  \left( \begin{array}{ccccc}
0 & 0 & 0 & \cdots & \chi(\varpi)\\
1 & 0 & 0 & \cdots & 0\\
0 & 1 & 0 & \cdots & 0\\
\vdots & \vdots & \vdots & \vdots & \vdots\\
0 & 0 & \cdots & 1 & 0
\end{array} \right)^{-1}
\left( \begin{array}{ccccc}
\chi(\varpi) & 0 & 0 & \cdots & 0\\
0 & 1 & 0 & \cdots & 0\\
0 & 0 & 1 & \cdots & 0\\
\vdots & \vdots & \vdots & \vdots & \vdots\\
0 & 0 & \cdots & 0 & 1
\end{array} \right) =
\left( \begin{array}{ccccc}
0 & 1 & 0 & \cdots & 0\\
0 & 0 & 1 & \cdots & 0\\
0 & 0 & 0 & \ddots & 0\\
\vdots & \vdots & \vdots & \vdots & 1\\
1 & 0 & \cdots & 0 & 0
\end{array} \right),  \]

which is an element of $SL(\ell,\mathbb{C})$.  Therefore, we may set

\[ \tau = \left( \begin{array}{ccccc}
\chi(\varpi) & 0 & 0 & \cdots & 0\\
0 & 1 & 0 & \cdots & 0\\
0 & 0 & 1 & \cdots & 0\\
\vdots & \vdots & \vdots & \vdots & \vdots\\
0 & 0 & \cdots & 0 & 1
\end{array} \right) \]

Then $\chi_{\tau} : X^{\sigma} \rightarrow \mathbb{C}^*$ is given by \[ \chi_{\tau}(\lambda_{(k,k, \cdots, k)}) = \lambda_{(k,k, \cdots, k)}(\tau) = (\chi(\varpi))^k \]
\qed

Recall that we have computed $\chi_{\phi}$ on $\mathfrak{o}_E^*$.  It remains to compute $\chi_{\phi}(\varpi)$.  Because of the isomorphism
$${}^0 T^{\Phi_{\sigma}} \times X^{\sigma} \stackrel{\sim}{\rightarrow} T^{\Phi_{\sigma}}$$ $$(\gamma, \lambda) \mapsto \gamma \lambda(\varpi),$$ we need to compute $\chi_{\phi}(1, \lambda_{(1,1, \cdots, 1)})$.

\begin{proposition}\label{kakapantspants}
Let $\ell = 2$.  Then $\chi_{\phi} = \chi \Delta_{\chi}$, where $\phi = Ind_{W_E}^{W_F}(\chi)$.
\end{proposition}

\proof
We have that $\chi_{\phi}(1, \lambda_{(1,1)}) = \chi_s(1) \chi_{\tau}(\lambda_{(1,1)}) = -\chi(\varpi)$.  Therefore, $\chi_{\phi}(\varpi) = -\chi(\varpi)$.  Recall that we have shown that $\chi_{\phi}|_{\mathfrak{o}_E^*} = \chi|_{\mathfrak{o}_E^*}$.  Since $\ell = 2$, $\Delta_{\chi}$ is the unique quadratic unramified character of $E^*$. Therefore, we have that $\Delta_{\chi}(\varpi) = -1$ and $\Delta_{\chi}|_{\mathfrak{o}_E^*} \equiv 1$, so $\chi_{\phi} = \chi \Delta_{\chi}$.
\qed

\begin{proposition}\label{kakapantspants2}
Let $\ell$ be an odd prime.  $\chi_{\phi} = \chi \Delta_{\chi}$.
\end{proposition}

\proof
We have that $\chi_{\phi}(1, \lambda_{(1,1, \cdots, 1)}) = \chi_s(1) \chi_{\tau}(\lambda_{(1,1, \cdots, 1)}) = \chi(\varpi)$.  Therefore, $\chi_{\phi}(\varpi) = \chi(\varpi)$.  Recall that we have shown that $\chi_{\phi}|_{\mathfrak{o}_E^*} = \chi|_{\mathfrak{o}_E^*}$.  Therefore,  $\chi_{\phi} = \chi$.  But recall that $\Delta_{\chi}$ is trivial since $\ell$ is an odd prime, so we have $\chi_{\phi} = \chi \Delta_{\chi}$.
\qed

\section{From a character of a torus to a representation for $GL(\ell,F)$}\label{glellfcharactertorus}

In this section we determine the representation that DeBacker--Reeder assign to a TRSELP for $GL(\ell,F)$, using the results from Section \ref{glellf}.  Note that $[X / (1 - w \theta) X]_{\mathrm{tor}} = 0$, so we may let $\lambda = 0$ (recall that $\lambda \in X_w$).  The proof of \cite[Lemma 2.7.2]{debackerreeder} implies that we may take $u_{\lambda} = 1$, and therefore $\Phi_{\lambda} = \Phi$.  It is also easy to see that we may take $w_{\lambda} = w$ (see \cite[section 2.7]{debackerreeder}) and $\dot{w}_{\lambda} = \dot{w}$, where $\dot{w}$ is a fixed choice of lift of $w$.  Since the theory of \cite{debackerreeder} is independent of any choices, we are free to choose a specific lift $\dot{w}$, which we do now.

Let $f(x)$ be a monic irreducible polynomial of degree $\ell$ over $\mathfrak{f}$.  Let $\tilde{f}(x)$ be a monic lift of $f(x)$ to $F[x]$.  We may write $E = F(\delta)$, where $\delta$ is a root of $\tilde{f}(x)$.   First set

\[ \tilde{w} := \left( \begin{array}{ccccc}
0 & 0 & 0 & \cdots & 1\\
1 & 0 & 0 & \cdots & 0\\
0 & 1 & 0 & \cdots & 0\\
\vdots & \vdots & \vdots & \vdots & \vdots\\
0 & 0 & \cdots & 1 & 0
\end{array} \right) \]

Recall that we need to find $p_{\lambda} \in G_{\lambda}$ such that $p_{\lambda}^{-1} \Phi(p_{\lambda}) = \dot{w}_{\lambda}$.  By choosing the basis $1, \delta, \delta^2, \cdots, \delta^{\ell-1}$ for $E$ over $F$, we may embed $E^*$ into $GL(\ell,F)$ in the standard way.  Denote this embedding $\varphi : E^* \hookrightarrow GL(\ell,F)$.

\begin{lemma}\label{conjugatingmatrix}
There exists an $A \in G_{\lambda}$ such that

\[A \left( \begin{array}{cccccc}
t & 0 & 0 & 0 & \cdots & 0\\
0 & \xi(t) & 0 & 0 & \cdots & 0\\
0 & 0 & \xi^2(t) & 0 & \cdots & 0\\
0 & 0 & 0 & \xi^3(t) & \cdots & 0\\
\vdots & \vdots & \vdots & \vdots & \ddots & \vdots\\
0 & 0 & \cdots & 0 & 0 & \xi^{\ell-1}(t)
\end{array} \right) A^{-1} = \varphi(t) \]

\noindent for all $t = a_0 + a_1 \delta + a_2 \delta^2 + \cdots + a_{\ell-1} \delta^{\ell-1} \in E^*$, $a_i \in F$.
\end{lemma}

\proof
Suppose $R(x)$ is a polynomial of degree $\ell$ in $F[x]$ that splits over $E$.  Then we get an isomorphism $$E[x]/(R(x)) \xrightarrow{\sim} \oplus_{i=1}^{\ell} E$$ $$p(x) \mapsto (p(a_1),p(a_2),...,p(a_{\ell}))$$ where $a_i$ are the roots of $R(x)$. Setting $R(x)$ to now be the minimal polynomial of $\delta$, and considering the basis $1, x, x^2, ..., x^{\ell-1}$ of $E[x] / (R(x))$ over $E$, we get an isomorphism $$E[x] / (R(x)) \xrightarrow{G} E \oplus E \oplus ... \oplus E$$ $$1 \mapsto (1,1,1,...,1)$$ $$x \mapsto (\delta, \Phi(\delta), \Phi^2(\delta), ..., \Phi^{\ell-1}(\delta))$$ $$x^2 \mapsto (\delta^2, \Phi(\delta)^2, \Phi^2(\delta)^2, ..., \Phi^{\ell-1}(\delta)^2)$$ $$...$$ $$x^{\ell-1} \mapsto (\delta^{\ell-1}, \Phi(\delta)^{\ell-1}, \Phi^2(\delta)^{\ell-1}, ..., \Phi^{\ell-1}(\delta)^{\ell-1})$$  This transformation yields the matrix

\[V := \left( \begin{array}{cccccc}
1 & \delta & \delta^2 & \delta^3 & \cdots & \delta^{\ell-1}\\
1 & \Phi(\delta) & \Phi(\delta)^2 & \Phi(\delta)^3 & \cdots & \Phi(\delta)^{\ell-1}\\
1 & \Phi^2(\delta) & \Phi^2(\delta)^2 & \Phi^2(\delta)^3 & \cdots & \Phi^2(\delta)^{\ell-1}\\
1 & \Phi^3(\delta) & \Phi^3(\delta)^2 & \Phi^3(\delta)^3 & \cdots & \Phi^3(\delta)^{\ell-1}\\
\vdots & \vdots & \vdots & \vdots & \vdots & \vdots\\
1 & \Phi^{\ell-1}(\delta) & \Phi^{\ell-1}(\delta)^2 & \Phi^{\ell-1}(\delta)^3 & \cdots &\Phi^{\ell-1}(\delta)^{\ell-1}
\end{array} \right) \]

\noindent We then set $A := V^{-1}$.  Note that what we have really done here is the following.  We first have taken the standard $E$-basis $e_1', e_2', ..., e_{\ell}'$ of $E \oplus E \oplus ... \oplus E$ and pulled it back by $G$ to get a basis $e_1, e_2, ..., e_{\ell}$ of $E[x] / (R(x))$. We have then shown that the standard embedding of an element $w \in E^*$ in $GL(\ell,F)$ (by considering its action on the basis $1, \delta, \delta^2, ..., \delta^{\ell-1}$ ) can be diagonalized over $E$ with respect to the basis $e_1, e_2, ..., e_{\ell}$.

Note that $V$ is a Vandermonde matrix.  Therefore, its determinant is $$\prod_{0 \leq i < j \leq \ell-1}  (\Phi^j(\delta) - \Phi^i(\delta)),$$ which has valuation zero.  Since we also have that the entries of $V$ are contained in $\mathfrak{o}_E^*$, we conclude that $V$, and hence $A$, is contained in $G_{\lambda}$.
\qed

Set $\tilde{s} = \tilde{w}^{-1} A^{-1} \Phi(A)$.  We now fix our choice of lift $\dot{w}$ of $w$ by setting $\dot{w} := \tilde{s}^{-1} \tilde{w} \tilde{s} \Phi(\tilde{s}) $, which we shall show is a legitimate lift.  We claim first that we may set $p_{\lambda} = A \tilde{s}$, and that $\tilde{s} \in G_{\lambda} \cap T$.  Note that since we will show that $\tilde{s} \in G_{\lambda} \cap T$, this shows that $p_{\lambda} \in G_{\lambda}$, which is required.  To prove all of this, consider the adjoint action of $A^{-1} \Phi(A)$ on $T$.  First, for $s \in T^{\Phi_w}$, we have $$(A^{-1} \Phi(A)) \cdot \Phi(s) = A^{-1} \Phi(A) \Phi(s) \Phi(A)^{-1} A = A^{-1} \Phi(A s A^{-1}) A = A^{-1} A s A^{-1} A = s$$
\noindent since Lemma \ref{conjugatingmatrix} implies that $AsA^{-1}$ is fixed by $\Phi$.

We therefore have that $(A^{-1} \Phi(A)) \cdot \Phi(s) = w \cdot \Phi(s) \ \forall s \in T^{\Phi_w}$.  Since $T^{\Phi_w}$ is dense in $T$ in the Zariski topology, we have that $(A^{-1} \Phi(A)) \cdot \Phi(s) = w \cdot \Phi(s) \ \forall s \in T$.  This implies that $( \tilde{w}^{-1} A^{-1} \Phi(A)) \cdot s = s \ \forall s \in T$ since $\tilde{w}$ is clearly a lift of $w$, which means that $$\tilde{w}^{-1} A^{-1} \Phi(A) s (\tilde{w}^{-1} A^{-1} \Phi(A))^{-1} = s \ \forall s \in T.$$ This means that $ \tilde{w}^{-1} A^{-1} \Phi(A) \in C_G(T) = T$, so in particular $ \tilde{w}^{-1} A^{-1} \Phi(A) = \tilde{s} \in T$.  But $A, \tilde{w} \in G_{\lambda}$ implies that $\tilde{w}^{-1} A^{-1} \Phi(A) \in G_{\lambda}$, which implies that $\tilde{s} \in G_{\lambda} \cap T$.  This shows that $p_{\lambda} \in G_{\lambda}$, which is required. Moreover, $p_{\lambda}^{-1} \Phi(p_{\lambda}) = (A \tilde{s})^{-1} \Phi(A \tilde{s}) = \tilde{s}^{-1} A^{-1} \Phi(A) \Phi(\tilde{s}) = \tilde{s}^{-1} \tilde{w} \tilde{s} \Phi(\tilde{s}) = \dot{w}$.  Finally, $\dot{w}$ is a lift of $w$ since $\tilde{w}$ is, and since $\tilde{s} \in T$, proving the claim.

Thus, we have a $p_{\lambda}$ such that $p_{\lambda}^{-1} \Phi_{\lambda}(p_{\lambda}) = \dot{w}$, and $\dot{w}$ is indeed a lift of $w$. Then if we define $T_{\lambda} := Ad(p_{\lambda}) T$, we get that $T_{\lambda}^{\Phi_{\lambda}}$ is the image of $E^*$ under $\varphi$.  This is crucial, since the depth zero supercuspidals of $GL(\ell,F)$ are constructed in section \ref{depthzeroclassicalconstruction} by first fixing an the embedding of $E^*$ into $GL(\ell,F)$.  The overall construction does not depend on the choice of embedding.  We have fixed the embedding $\varphi$.  DeBacker--Reeder are attaching a depth zero supercuspidal representation of $GL(\ell,F)$ to a Langlands parameter, and we need to show that their depth zero supercuspidal matches the depth zero supercuspidal attached in Theorem \ref{tamellc} (the latter of which, again, uses the construction in section \ref{depthzeroclassicalconstruction}, which assumes a fixed embedding, which we are assuming without loss of generality is $\varphi$).

Note that we have a simple description for the map $Ad(p_{\lambda})^{-1}$.
$$T_{\lambda}^{\Phi_{\lambda}} \xrightarrow{Ad(p_{\lambda})^{-1}} T^{\Phi_w}$$
$$\phi(t) \mapsto diag(t, \xi(t), \xi^2(t), ..., \xi^{\ell-1}(t)),$$ where $diag(d_1, d_2, ..., d_{\ell})$ denotes the diagonal $\ell$ by $\ell$ matrix with $d_1, d_2, ..., d_{\ell}$ on the diagonal, and where $t = a_0 + a_1 \delta + a_2 \delta^2 + \cdots + a_{\ell-1} \delta^{\ell-1}$.  Note that

$$T^{\Phi_w} = \{diag(a_0, \xi(a_0), \xi^2(a_0), ..., \xi^{\ell-1}(a_0)) : a_0 \in E^* \}$$
Finally, because of Propositions \ref{kakapantspants} and \ref{kakapantspants2} and the definition of $\chi_{\lambda}$, we have the following proposition.

\begin{proposition}
$\chi_{\lambda}$ is given by
\[\chi_{\lambda} (\varphi(t)) = \chi(t) \Delta_{\chi}(t) \]
for all $t = a_0 + a_1 \delta + a_2 \delta^2 + ... + a_{\ell-1} \delta^{\ell-1} \in E^*$.
\end{proposition}

Let us sum up the data that we have obtained so far.  Given a TRSELP for $GL(\ell,F)$, we have obtained a torus $T^{\Phi_w}$.  Given $\lambda = 0 \in X_w$, we have constructed $T_{\lambda}^{\Phi_{\lambda}}$ and $p_{\lambda}$.  We have $T_{\lambda}^{\Phi_{\lambda}} \cong E^*$.  From $\phi$ we have constructed a character $\chi_{\phi}$ of $T^{\Phi_w}$.  Via $Ad(p_{\lambda})$, we transported $\chi_{\phi}$ to a character $\chi_{\lambda}$ of $T_{\lambda}^{\Phi_{\lambda}}$.  We have shown that $\chi_{\lambda} = \chi \Delta_{\chi}$.  note that the restriction of $\chi_{\lambda}$ to ${}^0 T_{\lambda}^{\Phi_{\lambda}}$ factors through a character $\chi_{\lambda}^0$ of $\mathbb{T}_{\lambda}^{\Phi_{\lambda}}$.  Then, the packet of representations that DeBacker--Reeder construct in \cite{debackerreeder} from the data that we have obtained thus far is the single representation
$$Ind_{F^* GL(\ell,\mathfrak{o}_F)}^{GL(\ell,F)} (\chi_{\lambda} \otimes \kappa_{\lambda}^0) = \pi_{\chi \Delta_{\chi}}$$

\noindent Recall that in Section \ref{weilparameters}, the local Langlands correspondence for $GL(\ell,F)$, where $\ell$ is prime, was given as $$Ind_{W_E}^{W_F}(\chi) \mapsto \pi_{\chi \Delta_{\chi}}$$
We have therefore shown that the correspondence of DeBacker--Reeder coincides with the local Langlands correspondence.

\section{The positive depth correspondence of Reeder for $GL(\ell,F)$, $\ell$ an arbitrary prime}\label{positivedepthreederchapter}

In this section, we prove that the correspondence of \cite{reeder} agrees with the local Langlands correspondence of \cite{moy} for $GL(\ell,F)$, where $\ell$ is an arbitrary prime, if one assumes a certain compatibility condition, which we describe now.  Reeder's construction in \cite{reeder} begins by canonically attaching a certain admissible pair $(L/F, \Omega)$ to a Langlands parameter for $GL(\ell,F)$.  His construction then inputs this admissible pair into the theory of \cite{adler} in order to construct a supercuspidal representation $\pi(L, \Omega)$ of $GL(\ell,F)$.  The compatibility condition that we will need to assume is that $\pi(L, \Omega)$ is the same supercuspidal representation that is attached to $(L/F, \Omega)$ via the Howe construction in \cite{howe}.

Most of the arguments and setup are the same as in the depth-zero case, so there is not much to prove here.  We first very briefly review the construction of \cite{reeder} and refer to \cite{reeder} for various definitions and notions that are not explained here.

\subsection{Review of construction of Reeder}\label{reederconstruction}

Let $\mathbf{G}$ be an $F$-quasisplit and $F^u$-split connected reductive group.  Let $\mathbf{B} \subset \mathbf{G}$ be a Borel subgroup defined over $F$, and $\mathbf{T}$ a maximal torus of $\mathbf{B}$.

The Langlands parameters considered in \cite{reeder} are the maps $\phi : W_F \rightarrow {}^L G = \ <\hat{\theta}> \ltimes \hat{G}$ such that

(1) $\phi$ is trivial on $I^{(r+1)}$ and nontrivial on $I^{(r)}$ for some integer $r > 0$.  Here, $\{ I^{(k)} \}_{k \geq 0}$ is a filtration on $I_F$ defined in \cite[section 5.2]{reeder}.

(2) The centralizer of $\phi(I^{(r)})$ in $\hat{G}$ is a maximal torus of $\hat{G}$.  This condition is the \emph{regularity} condition.

(3) $\phi(\Phi) \in \hat{\theta} \ltimes \hat{G}$, and the centralizer of $\phi(W_F)$ in $\hat{G}$ is finite, modulo $\hat{Z}^{\hat{\theta}}$.  This is the \emph{ellipticity} condition.

We may conjugate $\phi$ by an element of $\hat{G}$ so that $\phi(I_F) \subset \hat{T}$, and $\phi(\Phi) = \hat{\theta} f$, where $f \in \hat{N}$.  Let $\hat{w}$ be the image of $f$ in $\hat{W}_o$, and let $w$ be the element of $W_o$ dual to $\hat{w}$.  We say that the element $w$ is \emph{associated} to $\phi$.

Set $\sigma = w \theta$ and suppose its action on $X$ has order $n$.  From an above such Langlands parameter, Reeder defines a $\hat{T}$-conjugacy class of Langlands parameters $$\phi_T : W_F \rightarrow {}^L T_{\sigma}$$ in the exact same way as in the depth-zero case.  In particular, the element $\tau$ is defined in the same way.

As in the depth-zero case, a bijection is later given between $\hat{T}$-conjugacy classes of continuous homomorphisms $$\phi : W_F / I^{(r+1)} \rightarrow {}^L T_{\sigma}$$ for which $\phi(\Phi) \in \hat{\sigma} \ltimes \hat{T}$ and characters of $T^{\Phi_{\sigma}}$ that are trivial on $T_{r+1}^{\Phi_{\sigma}}$, where $\{ T_k \}_{k \geq 0}$ is the canonical filtration on $T$ (see \cite[section 5.3]{reeder}).  This is done as follows.  We have a composite isomorphism (see \cite[section 5.3]{reeder})

\begin{equation*}
\mathrm{Hom}_{Ad(\Phi), \hat{\sigma}}(I_F / I^{(r+1)}, \hat{T}) \cong \mathrm{Hom}_{Ad(\Phi), \hat{\sigma}}(I_F / I_n^{(r+1)}, \hat{T}) = \mathrm{Hom}_{\Phi, \hat{\sigma}}(\mathfrak{o}_n^* / (1 + \mathfrak{p}_n^{r+1}), \hat{T}) =
\end{equation*}
\begin{equation}\label{compositeisomorphismreedercool}
\mathrm{Hom}_{\Phi_{\sigma}, Id}(X \otimes (\mathfrak{o}_n^* / (1 + \mathfrak{p}_n^{r+1})), \mathbb{C}^*) = \mathrm{Hom}_{\Phi_{\sigma}, Id}({}^0 T^{\Phi_{\sigma}^n} / T_{r+1}^{\Phi_{\sigma}^n}, \mathbb{C}^*) = \mathrm{Hom}({}^0 T^{\Phi_{\sigma}} / T_{r+1}^{\Phi_{\sigma}}, \mathbb{C}^*)
\end{equation}

Under this composite isomorphism, $s := \phi|_{I_F}$ maps to a character $\chi_s \in \mathrm{Hom}({}^0 T^{\Phi_{\sigma}} / T_{r+1}^{\Phi_{\sigma}}, \mathbb{C}^*)$.  Then, if $\phi(\Phi) = \hat{\sigma} \ltimes \tau$, we get that $\tau$ gives rise to a character of $X^{\sigma}$ given by $\chi_{\tau}(\lambda) := \lambda(\tau)$ for $\lambda \in X^{\sigma}$, just as in the depth-zero case.  Recalling that $T^{\Phi_{\sigma}} = {}^0 T^{\Phi_{\sigma}} \times X^{\sigma}$, we define a character $\chi_{\phi}$ of $T^{\Phi_{\sigma}}$ by $\chi_{\phi} := \chi_s \otimes \chi_{\tau}$, which is our desired character of $T^{\Phi_{\sigma}}$ constructed from the Langlands parameter $\phi$.

As in the depth-zero case, we have the set $X_w$.  To $\lambda \in X_w$, Reeder associates a 1-cocycle $u_{\lambda}$, hence a twisted Frobenius $\Phi_{\lambda} = Ad(u_{\lambda}) \circ \Phi$.  Moreover, to $\lambda$ is associated an affine Weyl group element $w_{\lambda}$, a parahoric subgroup $G_{x_{\lambda}}$, and an element $p_{\lambda} \in G_{x_{\lambda}}$ such that $p_{\lambda}^{-1} \Phi_{\lambda}(p_{\lambda})$ is a lift of $w_{\lambda}$.  We then define $T_{\lambda} := Ad(p_{\lambda}) T$ and set $\chi_{\lambda} := \chi_{\phi} \circ Ad(p_{\lambda})^{-1}$.  To the torus $T_{\lambda}$ and the character $\chi_{\lambda}$, we apply the construction of \cite{adler} to obtain a supercuspidal representation.  Then, Reeder constructs a packet $\Pi(\phi)$ of representations on the pure inner forms of $G$, parameterized by $\mathrm{Irr}(C_{\phi})$, using the above construction.

\subsection{The case of $GL(\ell,F)$, where $\ell$ is prime}

In this section, we consider the group $\mathbf{G}(F) = GL(\ell,F)$, where $\ell$ prime.  Let $\phi : W_F \rightarrow {}^L G$ be one of the Langlands parameters for $\mathbf{G}(F) = GL(\ell,F)$ that is considered in section \ref{reederconstruction}.

\begin{lemma}
$\phi = Ind_{W_E}^{W_F}(\chi)$ for some admissible pair $(E/F, \chi)$, where $\chi$ has positive level and $E/F$ is degree $\ell$ unramified.
\end{lemma}

\proof
The proof is similar as in the $GL(2,F)$ case, but we include it for completeness purposes.  As in the depth-zero case in Section \ref{glellf}, we may conjugate $\phi$ by an element of $\hat{G}$ so that the Weyl group element $w$ that is associated to $\phi$ is the Weyl group element $(1 \ 2 \ 3 \ ... \ \ell)$ in the symmetric group on $\ell$ letters.  We know that $\phi$ is an irreducible admissible $\phi : W_F \rightarrow GL(\ell,\mathbb{C})$ that is trivial on $I^{(r+1)}$ and nontrivial on $I^{(r)}$ for some integer $r > 0$.  Let $E$ be the degree $\ell$ unramified extension of $F$.  Again, any representation $Ind_{W_E}^{W_F}(\Omega)$ where $(E/F, \Omega)$ is an admissible pair is equivalent to the representation $\kappa : W_F \rightarrow GL(\ell,\mathbb{C})$ satisfying:

1) $\kappa|_{W_E}$ is given by $\Omega \in \widehat{E^*}$ by the local Langlands correspondence for tori.

2) \[ \kappa(\Phi) = \left( \begin{array}{ccccc}
0 & 0 & 0 & \cdots & \Omega(\varpi)\\
1 & 0 & 0 & \cdots & 0\\
0 & 1 & 0 & \cdots & 0\\
\vdots & \vdots & \ddots & \vdots & \vdots\\
0 & 0 & \cdots & 1 & 0
\end{array} \right) \]

We want to show that $\phi$ satisfies the two conditions above, for some admissible pair $(E/F, \chi)$.  Let's restrict $\phi$ to $W_E$.  By the composite isomorphism (\ref{compositeisomorphismreedercool}), $\phi|_{I_E}$ gives rise to a character $\ddot{\chi}$ of $\mathfrak{o}_E^*$.  Then, by following the composite isomorphism (\ref{compositeisomorphismreedercool}) backwards, one sees that

\[\phi(x) = \left( \begin{array}{ccccc}
\ddot{\chi}(r_{\ell}(x)) & 0 & 0 & \cdots & 0\\
0 & \ddot{\chi}^{\xi}(r_{\ell}(x)) & 0 & \cdots & 0\\
0 & 0 & \ddot{\chi}^{\xi^2}(r_{\ell}(x)) & \cdots & 0\\
\vdots & \vdots & \vdots & \ddots & \vdots\\
0 & 0 & 0 & \cdots & \ddot{\chi}^{\xi^{\ell-1}}(r_{\ell}(x))
\end{array} \right) \]

\noindent as in the depth-zero case.  Now, as in Propositions \ref{chitaugl2} and \ref{chitauglelldepthzero}, we know that

\[ \phi(\Phi) = \left( \begin{array}{ccccc}
0 & 0 & 0 & \cdots & a\\
1 & 0 & 0 & \cdots & 0\\
0 & 1 & 0 & \cdots & 0\\
\vdots & \vdots & \ddots & \vdots & \vdots\\
0 & 0 & \cdots & 1 & 0
\end{array} \right) \]

\noindent for some $a \in \mathbb{C}^*$, because of the ellipticity condition on $\phi$.  Therefore, we have that

\[ \phi(\Phi_E) = \phi(\Phi^{\ell}) = \phi(\Phi)^{\ell} = \left( \begin{array}{ccccc}
a & 0 & 0 & \cdots & 0\\
0 & a & 0 & \cdots & 0\\
0 & 0 & a & \cdots & 0\\
\vdots & \vdots & \vdots & \ddots & \vdots\\
0 & 0 & \cdots & 0 & a
\end{array} \right). \]

\noindent  Then $\ddot{\chi}$ extends to a character, denoted $\chi$, of $E^*$, by setting $\chi(\varpi) := a$ and $\chi|_{\mathfrak{o}_E^*} := \ddot{\chi}|_{\mathfrak{o}_E^*}$. One can now see that $\phi = Ind_{W_E}^{W_F}(\chi)$.  By the regularity condition on $\phi$, we get that $\ddot{\chi} \neq \ddot{\chi}^{\xi}$, and thus $(E/F, \chi)$ is an admissible pair.  Finally, $\chi$ has positive level since $r > 0$.
\qed

\begin{proposition}\label{positivedepthcharactergl2}
Let $\ell = 2$.  Then $\chi_{\phi} = \chi \Delta_{\chi}$.
\end{proposition}

\proof
The analogous arguments as in the depth-zero case show that $\chi_{\phi}|_{\mathfrak{o}_E^*} = \chi|_{\mathfrak{o}_E^*}$.  In particular, let $z \in \mathfrak{o}_E^*$.  Let $x \in I_F$ be any element such that $r_2(x) = z$ (where $r_2$ is as in \cite[Section 5.1]{reeder}), and let $\Gamma$ be the cocharacter $$t \mapsto \mat{t}{0}{0}{1}.$$  Then $$N_{\sigma}(\Gamma \otimes r_2(x)) = \mat{z}{0}{0}{\overline{z}}.$$ Moreover, by the same arguments as in Proposition \ref{whatevers2}, we get $\phi(x) = \mat{\chi \left(r_2(x) \right)}{0}{0}{\chi \left(\overline{r_2(x)} \right)}$ so that $\Gamma(\phi(x)) = \chi(z)$, where here we are viewing $\Gamma$ as a character of $\hat{T}$.  Finally, as we may take $\tau$ to be the same element as in the depth-zero case, we have that $\chi_{\phi}(\varpi) = -\chi(\varpi)$, so that $\chi_{\phi} = \chi \Delta_{\chi}$.
\qed

\begin{proposition}\label{glellcharactertorusreeder}
Let $\ell$ be an odd prime.  Then $\chi_{\phi} = \chi \Delta_{\chi}$.
\end{proposition}

\proof
The analogous proof as in Proposition \ref{positivedepthcharactergl2} and the depth-zero case works here.
\qed

Note that $[X / (1 - w \theta) X]_{\mathrm{tor}} = 0$, so we may let $\lambda = 0$ (recall that $\lambda \in X_w$).  It is easy to see that we may again take $u_{\lambda} = 1$, and therefore $\Phi_{\lambda} = \Phi$.  It is also easy to see that we may take $w_{\lambda} = w$ (see \cite[section 6.4]{reeder}), and we may also take the same $p_{\lambda}$ as in the depth-zero case in Section \ref{glellfcharactertorus}.  So we have the same $T_{\lambda}$ as in section \ref{glellfcharactertorus} and the analogous $\chi_{\lambda}$.

We have therefore shown that if we assume the compatibility condition in the beginning of section \ref{positivedepthreederchapter}, then by Proposition \ref{glellcharactertorusreeder}, the Reeder construction attaches the representation $\pi_{\chi \Delta_{\chi}}$ to the Langlands parameter $\phi = Ind_{W_E}^{W_F}(\chi)$. This shows that as long as we assume this compatibility condition, the correspondences of \cite{reeder} and \cite{moy} agree for $GL(\ell,F)$, where $\ell$ is an odd prime.

Acknowledgements: This paper has benefited from conversations with Jeffrey Adams, Jeffrey Adler, Gordan Savin, and Jiu-Kang Yu.


\begin{thebibliography}{9}

\bibitem{adler}
  Jeffrey Adler,
  \emph{Refined anisotropic K-types and supercuspidal representations}, Pacific J. Math., 185 (1998), pp. 1-32.

\bibitem{bushnellhenniart}
  Colin Bushnell and Guy Henniart,
  \emph{The Local Langlands Conjecture for GL(2)}, A Series of Comprehensive Studies in Mathematics, Volume 335, Springer Berlin Heidelberg, 2006.

\bibitem{bushnellhenniart1}
  Colin Bushnell and Guy Henniart,
  \emph{The Essentially Tame Local Langlands Correspondence, I},
  J. Amer. Math. Soc. 18 (2005), no. 3, 685--710.

\bibitem{carter}
  Roger Carter,
  \emph{Finite Groups of Lie Type: Conjugacy Classes and Complex Characters}, John Wiley and Sons Inc, 1993.

\bibitem{debacker}
  Stephen DeBacker,
  \emph{On Supercuspidal Characters of $GL_{\ell}$, $\ell$ a prime}, Ph.D. thesis, University of Chicago, 1997.

\bibitem{debackerreeder}
  Stephen DeBacker and Mark Reeder,
  \emph{Depth-zero supercuspidal $L$-packets and their stability.}
  Ann. of Math. (2) 169 (2009), no. 3, 795--901.

\bibitem{dignemichel}
  Francois Digne, Jean Michel,
  \emph{Representations of Finite Groups of Lie Type}, London Mathematical Society Student Texts 21, Cambridge University Press, 1991.

\bibitem{howe}
  Roger Howe,
  \emph{Tamely ramified supercuspidal representations of $GL_n(F)$},
   Pacific Journal of Math.  73  (1977),  437--460.

\bibitem{moy}
  Allen Moy,
  \emph{Local Constants and the Tame Langlands Correspondence},
   American Journal of Math.  108  (1986),  no. 4, 863--929.

\bibitem{reeder}
  Mark Reeder,
  \emph{Supercuspidal $L$-packets of positive depth and twisted Coxeter elements.} J. Reine Angew. Math. 620 (2008), 1-33.

\bibitem{serre}
  Jean-Pierre Serre,
  \emph{A Course in Arithmetic}, 4th edition, Graduate Texts in Mathematics, Springer-Verlag Berlin and Heidelberg GmbH and Co. K, 1996.

\bibitem{serre1}
  Jean-Pierre Serre,
  \emph{Local Fields}, Graduate Texts in Mathematics, Springer-Verlag New York, 1979.

\bibitem{tate}
  John Tate,
  \emph{Number theoretic background.} Automorphic forms, representations and L-functions (Proc. Sympos. Pure Math., Oregon State Univ., Corvallis, Ore., 1977), Part 2, pp. 3-26, Proc. Sympos. Pure Math., XXXIII, Amer.
Math. Soc., Providence, RI, 1979.

\end{thebibliography}
\end{document}